\documentclass{emsprocart}
\usepackage{hyperref}

\newcommand{\mlabel}[1]{\marginpar{#1}\label{#1}}
\renewcommand{\mlabel}{\label}

\renewcommand{\:}{\colon}

\newcommand{\la}{\langle}
\newcommand{\ra}{\rangle}
\newcommand{\oline}{\overline}
\newcommand{\subeq}{\subseteq} 
 
\newcommand{\into}{\hookrightarrow}

\newcommand{\res}{\vert}
\newcommand{\Part}{\mathop{{\rm Part}}\nolimits}
\newcommand{\id}{\mathop{{\rm id}}\nolimits}
\newcommand{\supp}{\mathop{{\rm supp}}\nolimits}

\newcommand{\eps}{\varepsilon}
\newcommand{\onto}{\to\mskip-14mu\to} 

\newcommand{\1}{\mathbf{1}}
\newcommand{\bV}{{\mathbb V}}

\newcommand{\bx}{{\bf x}}
\newcommand{\bn}{{\bf n}}
\newcommand{\bd}{{\bf d}}
\renewcommand{\L}{\mathop{\bf L{}}\nolimits}
\newcommand{\derat}[1]{\frac{d}{dt} \hbox{\vrule width0.5pt
                height 5mm depth 3mm${{}\atop{{}\atop{\scriptstyle t=#1}}}$}}

\newcommand{\g}{\mathfrak g}
\newcommand{\fb}{\mathfrak b}
\newcommand{\fk}{\mathfrak k}
\newcommand{\fq}{\mathfrak q}
\newcommand{\fp}{\mathfrak p}
\newcommand{\ft}{\mathfrak t}
\newcommand{\fh}{\mathfrak h}
\newcommand{\fu}{\mathfrak u}

\newcommand{\fz}{\mathfrak z}
\newcommand{\z}{\mathfrak z}
\newcommand{\fsp}{\mathfrak {sp}} 
 
\newcommand{\gau}{\mathfrak {gau}} 
 
\newcommand{\vir}{\mathfrak {vir}} 

\newcommand{\C}{{\mathbb C}}
\newcommand{\T}{{\mathbb T}}
\newcommand{\R}{{\mathbb R}}
\newcommand{\Z}{{\mathbb Z}}
\newcommand{\N}{{\mathbb N}}
\newcommand{\bP}{{\mathbb P}}
\newcommand{\bS}{{\mathbb S}}

\newcommand{\K}{{\mathbb K}}

\newcommand{\cA}{\mathcal{A}} 
\newcommand{\cD}{\mathcal{D}} 
 
\newcommand{\cH}{\mathcal{H}}
\newcommand{\cK}{\mathcal{K}}
\newcommand{\cP}{\mathcal{P}}
\newcommand{\cL}{\mathcal{L}}
\newcommand{\cM}{\mathcal{M}}
\newcommand{\cV}{\mathcal{V}}
\newcommand{\cS}{\mathcal{S}}
\newcommand{\cW}{\mathcal{W}}

\newcommand{\ad}{\mathop{{\rm ad}}\nolimits}

\newcommand{\Ad}{\mathop{{\rm Ad}}\nolimits}
\newcommand{\conv}{\mathop{{\rm conv}}\nolimits}
\newcommand{\Aut}{\mathop{{\rm Aut}}\nolimits}
\newcommand{\Heis}{\mathop{{\rm Heis}}\nolimits}
\newcommand{\Diff}{\mathop{{\rm Diff}}\nolimits}

\newcommand{\Ext}{\mathop{{\rm Ext}}\nolimits}
\newcommand{\ev}{\mathop{{\rm ev}}\nolimits}
\newcommand{\GL}{\mathop{{\rm GL}}\nolimits}
\newcommand{\Gau}{\mathop{{\rm Gau}}\nolimits}

\newcommand{\SO}{\mathop{{\rm SO}}\nolimits}

\newcommand{\Sp}{\mathop{{\rm Sp}}\nolimits}
\newcommand{\Spann}{\mathop{{\rm span}}\nolimits}
\newcommand{\Hom}{\mathop{{\rm Hom}}\nolimits}

\renewcommand{\Im}{\mathop{{\rm Im}}\nolimits}

\newcommand{\Spec}{\mathop{{\rm Spec}}\nolimits}
\newcommand{\U}{\mathop{\rm U{}}\nolimits}

\newcommand{\OO}{\mathop{\rm O{}}\nolimits}
\newcommand{\Vir}{\mathop{\rm Vir{}}\nolimits}
\newcommand{\dd}{{\tt d}} 

\renewcommand{\tilde}{\widetilde}
\renewcommand{\hat}{\widehat}

\renewcommand{\phi}{\varphi} 

\newcommand{\pmat}[1]{\begin{pmatrix} #1 \end{pmatrix}}

\contact[neeb@math.fau.de]{Karl-Hermann Neeb, Department of Mathematics, 
Friedrich-Alexander-University, Erlangen-Nuremberg, 
Cauerstrasse 11, 91056 Erlangen, Germany} 

\newtheorem{theorem}{Theorem}[section]

\newtheorem*{coro}{Corollary} 

\theoremstyle{definition}
\newtheorem{definition}[theorem]{Definition}
\newtheorem{remark}[theorem]{Remark}
\newtheorem{example}[theorem]{Example} 

\newcommand{\Ker}{\operatorname{Ker}}

\title[Bounded and Semibounded Representations of Infinite Dimensional Lie Groups]
{Bounded and Semibounded Representations\\ of Infinite Dimensional Lie Groups}

\author[Karl-Hermann Neeb]{Karl-Hermann Neeb\thanks{The author gratefully 
acknowledges the support of DFG-grant NE 413/7-2 in the framework
of SPP ``Representation Theory''.}} 
\begin{document}

\begin{abstract} In this note we describe the recent progress in 
the classification of bounded and semibounded representations of 
infinite dimensional Lie groups. We start with a discussion 
of the semiboundedness condition and how the new concept 
of a smoothing operator can be used to construct $C^*$-algebras 
(so called host algebras) 
whose representations are in one-to-one correspondence with certain 
semibounded representations of an infinite dimensional
 Lie group $G$. This makes the full power of $C^*$-theory available 
in this context. 

Then we discuss the classification of bounded representations of 
several types of unitary groups on Hilbert spaces and of gauge groups. 
After explaining the method of holomorphic induction as a means to 
pass from bounded representations to semibounded ones, we describe 
the classification of semibounded representations for 
hermitian Lie groups of operators, 
loop groups (with infinite dimensional targets), the Virasoro group 
and certain infinite dimensional oscillator groups. 
\end{abstract}

\begin{classification}
Primary 22E65; Secondary 22E45
\end{classification}

\begin{keywords}
infinite dimensional Lie group, host algebra, 
semibounded representation, bounded representation, 
holomorphic induction 
\end{keywords}

\maketitle

\section{Introduction} 

This article surveys some of the recent developments in 
the representation theory of infinite dimensional 
Lie groups. Here an infinite dimensional Lie group is a group $G$ 
which is a smooth manifold modeled on a locally convex space for which 
the group operations are smooth. As usual, we  identify the 
Lie algebra $\g$ of $G$ with the tangent space $T_\1(G)$ in the neutral 
element~$\1$. We shall assume throughout that $G$ has an 
{\it exponential function}, i.e., a smooth map 
$\exp \: \g \to G$ such  that, for every $x \in \g$, 
the curve $\gamma_x(t) := \exp(tx)$ is a one-parameter group 
with $\gamma_x'(0) = x$. We refer to \cite{Ne06} for the 
basics of infinite dimensional Lie theory. 

Typical examples of infinite dimensional Lie groups are 
\begin{itemize}
\item {\it Banach--Lie groups}, such as the group $\GL(E)$ of invertible 
bounded linear operators on a Banach space $E$, or, more generally, 
the group $\cA^\times$ of units in a unital Banach algebra $\cA$. 
For $G = \cA^\times$, the Lie algebra $\g$ can be identified with $\cA$, endowed with 
the commutator bracket, and
 the exponential function is given by the convergent 
exponential series $\exp(x) = \sum_{n = 0}^\infty \frac{x^n}{n!}$. 
\item {\it Mapping groups}, such as the group $G = C^\infty_c(M,K)$ of compactly 
supported smooth functions on a finite dimensional $\sigma$-compact 
manifold with values in a Lie group~$K$. 
\begin{footnote}
{More general examples of this 
type arise as groups of compactly supported gauge transformations 
of $K$-principal bundles $P \to M$ (cf.\ Section~\ref{sec:6}).}
\end{footnote}
Then $\g = C^\infty_c(M,\fk)$, with the Lie bracket 
$[\xi,\eta](m) := [\xi(m),\eta(m)]$, 
and the exponential function is given by 
$\exp(\xi) := \exp_K \circ \xi$ for $\xi \in\g$. 
\item {\it Diffeomorphism groups}, such as the group 
$\Diff_c(M)$ of compactly supported diffeomorphisms of a finite 
dimensional manifold. Then the Lie algebra is the space $\cV_c(M)$ of 
compactly supported smooth vector fields and the exponential 
function is given by the time-$1$-flow map 
$\exp(X) = \Phi^X_1$, where 
$(t,m) \mapsto \Phi^X_t(m)$ is the flow of the vector field $X$. 
\end{itemize}

We are interested in {\it unitary representations}, i.e., 
homomorphisms \break $\pi \: G \to \U(\cH)$ into the unitary group 
of a complex Hilbert space $\cH$ for which all orbit maps 
$\pi^v \:  G \to \cH, g \mapsto \pi(g)v$ are continuous. 
The exponential function permits us  
to associate to every unitary representation $(\pi, \cH)$ 
and every $x \in \g$ the unitary one-parameter group  
$\pi_x(t) := \pi(\exp tx)$. By Stone's Theorem, there exists a 
selfadjoint operator $A := -i\oline{\dd\pi}(x)$ on $\cH$ with 
\begin{equation}
  \label{eq:stone}
\pi_x(t) = \pi(\exp tx) = e^{t \oline{\dd\pi}(x)} = e^{itA} \quad \mbox{ for } 
\quad t \in \R, 
\end{equation}
where the exponential on the right is to be understood in the sense of 
functional calculus of selfadjoint operators.

To develop a reasonably general theory of unitary representations 
of infinite dimensional Lie groups, new approaches have to be developed 
which neither rest on the fine structure theory available 
for finite dimensional groups nor on the existence of 
invariant measures. One thus has 
to specify suitable classes of representations for which 
it is possible to develop sufficiently powerful tools. 
The notion of semiboundedness specified below 
is very much inspired by the concepts and requirements of 
mathematical physics and provides a unifying framework 
for a substantial class of representations and several interesting phenomena. 

The selfadjoint operators $\oline{\dd\pi}(x)$ permit us to formulate 
suitable regularity conditions for unitary representations 
to specify interesting classes of representations for which a powerful  
theory can be developed. To ensure that Lie theoretic 
methods can be applied, we have to require that the representation 
$(\pi, \cH)$ is {\it smooth}, i.e., that the subspace 
\[ \cH^\infty  := \cH^\infty(\pi) := \{ v \in \cH \: \pi^v \in C^\infty(G,\cH) \} \] 
is dense in $\cH$.
\begin{footnote}
{If $G$ is finite dimensional, then 
L.~G\aa{}rding's Theorem asserts that this is always 
satisfied \cite{Ga47}, but this is not the case for infinite dimensional Lie 
groups. The representation of the additive 
Banach--Lie group $G := L^2([0,1],\R)$ on 
$\cH = L^2([0,1],\C)$ by $\pi(g)f := e^{ig}f$ is continuous 
with $\cH^\infty = \{0\}$ (\cite{BN08}). }  
\end{footnote}
For a smooth unitary representation $(\pi,\cH)$, we consider its 
{\it support functional} 
\[ s_\pi \: \g \to \R\cup \{\infty\}, \quad 
s_\pi(x) := \sup\big(\Spec(i\oline{\dd\pi}(x))\big).\] 
This is a lower semicontinuous (its epigraph is closed), 
positively homogeneous convex function 
on $\g$ which is invariant under the adjoint action of $G$ on $\g$ 
(cf.\ Section~\ref{sec:2}). Its natural ``domain of regularity'' is the open 
convex cone 
\[ W_\pi := \{ x \in \g \:(\exists U\ \text{open},x \in U)\  
\sup(s_\pi\res_U) < \infty \}.\] 
We say that $\pi$ is {\it semibounded} if $W_\pi$ is non-empty and 
that $\pi$ is {\it bounded} if $W_\pi = \g$ (cf.\ \cite{Ne09}).

It turns out that, to a large extent, the methods of classical harmonical 
analysis on locally compact groups can be extended to 
semibounded unitary representations of infinite dimensional Lie groups. 
Below we shall discuss some of the cornerstones of this theory: 
\begin{itemize}
\item The existence of $C^*$-algebras $\cA$ whose representations are 
in one-to-one correspondence with those representations of $G$ 
satisfying a certain spectral condition (Theorem~\ref{thm:5.2}). 
Here the method of smoothing 
operators developed recently in \cite{NSZ15} is the key to complete the picture 
(Section~\ref{sec:3}).  
\item The Spectral Theorem for semibounded representations of 
locally convex spaces, resp., abelian Lie groups  \cite[Thm.~4.1]{Ne09} 
(Remark~\ref{rem:4.4}).  
\item The method of holomorphic induction which is presently the most effective 
tool to classify semibounded representations 
(see \cite{Ne13} for the case of Banach--Lie groups and 
\cite[App.~C]{Ne14} for generalizations to Fr\'echet--Lie groups) 
(Section~\ref{sec:7}).   
\item Some classification results for bounded representations 
that are particularly important as input for holomorphic induction 
(\cite{Ne98, Ne14b, JN14, BN14}) (Sections~\ref{sec:5}, \ref{sec:6}). 
\item Some classification results for semibounded representations of 
various classes of Lie groups 
(\cite{Ne12, Ne14, NS14}) (Sections~\ref{sec:8}, \ref{sec:9}, \ref{sec:10}). 
\end{itemize}

In addition to their tractability, the 
restriction to the class of semibounded representations 
is to a large extent motivated by the fact that 
representations arising in quantum physics 
are often characterized by the 
requirement that the Lie algebra $\g$ 
of $G$ contains an element $h$, corresponding 
to the Hamiltonian of the underlying physical system, for which the 
spectrum of the operator $-i \dd \pi(h)$ 
is bounded from below. These representations are called {\it positive energy 
representations} (for the appearance of such conditions is physics, see 
\cite{Se67, Se78}, 
\cite{SeG81}, \cite{Ca83}, \cite{Mick89}, \cite{PS86}, 
\cite{CR87}, \cite{Bo96}, \cite{FH05}, \cite{Bak07}). 
Semiboundedness is a  stable version of the 
positive energy condition. It 
means that the selfadjoint operators $i\dd\pi(x)$ from the derived 
representation are uniformly bounded from below for all $x$ in some 
non-empty open subset of $\g$. 


\section{Momentum sets of smooth unitary representations} \mlabel{sec:2}

In this section we introduce the momentum map and the momentum 
set of a smooth unitary representation of a Lie group~$G$. 

Let $(\pi, \cH)$ be a smooth unitary representation of $G$, so that the 
subspace $\cH^\infty$ of smooth vectors is dense. On  $\cH^\infty$ the 
derived representation $\dd\pi$ of the Lie algebra $\g = \L(G)$ 
is defined by 
\[ \dd\pi(x)v := \derat0 \pi(\exp tx)v.  \] 
The invariance of $\cH^\infty$ under $\pi(G)$ implies that, 
for $x \in \g$, the operator $i \dd\pi(x)$ on $\cH^\infty$ 
is  essentially selfadjoint (cf.\ 
\cite[Thm.~VIII.10]{RS75}) 
and that its closures coincides with the selfadjoint generator 
$i\oline{\dd\pi}(x)$ of the unitary one-parameter group 
$\pi_x(t) := \pi(\exp tx)$. 

\begin{definition} (a) Let $\bP({\cal H}^\infty) = \{ [v] := 
\C v \: 0 \not= v \in 
{\cal H}^\infty\}$ 
denote the projective space of the subspace ${\cal H}^\infty$ 
of smooth vectors. The map 
\[  \Phi_\pi \: \bP({\cal H}^\infty)\to \g' \quad \hbox{ with } \quad 
\Phi_\pi([v])(x) 
= \frac{1}{i}  \frac{\la  \dd\pi(x)v, v \ra}{\la v, v \ra} \] 
is called the {\it momentum map of the unitary representation $\pi$}.  
The operator $i\dd\pi(x)$ is symmetric so
that the right hand side is real, and since $v$ is a smooth vector, 
it defines a continuous linear functional on $\g$. 
We also observe that we have a natural action of $G$ on 
$\bP(\cH^\infty)$ by $g.[v] := [\pi(g)v]$, and the relation 
$$ \pi(g) \dd\pi(x) \pi(g)^{-1} = \dd\pi(\Ad(g)x) $$
immediately implies that $\Phi_\pi$ is equivariant with respect 
to the coadjoint action of $G$ on the topological dual space $\g'$.
\begin{footnote}
{If $G$ is a Banach--Lie group, then $\cH^\infty$ carries a natural 
Fr\'echet structure for which the $G$-action on the complex 
Fr\'echet--K\"ahler manifold $\bP(\cH^\infty)$, endowed with the 
Fubini--Study metric, is Hamiltonian with momentum 
map $\Phi_\pi$ (cf.\ \cite[Thm.~4.5]{Ne10}). This motivates the 
terminology.}  
\end{footnote}

(b) The weak-$*$-closed convex hull 
$I_\pi \subeq \g'$ of the image of $\Phi_\pi$ is called the 
{\it (convex) momentum set of $\pi$}. In view of the equivariance 
of $\Phi_\pi$, it is an $\Ad^*(G)$-invariant closed convex subset of $\g'$. 
\end{definition}

The momentum set $I_\pi$ 
provides complete information on the extreme spectral 
values of the operators $i\dd\pi(x)$: 
\begin{equation}
  \label{eq:momspec}
\sup(\Spec(i\dd\pi(x))) = s_{\pi}(x)= \sup  \la I_\pi,- x \ra 
 \quad \mbox{ for } \quad x \in \g 
\end{equation}
(cf.\ \cite[Lemma 5.6]{Ne08}). 
This relation shows that $s_\pi$ is the {\it support functional} 
of the convex subset $I_\pi \subeq \g'$, which implies that 
it is lower semicontinuous and convex. It is obviously positively homogeneous. 
The semibounded unitary representations are those for which 
$s_\pi$, resp., the  momentum set $I_\pi$, contains significant 
information on~$\pi$.  For these the set $I_\pi$ 
is {\it  semi-equicontinuous} in the sense that its support function 
$s_\pi$ is bounded in a neighborhood of some $x_0 \in \g$.

\begin{remark}(Physical interpretation) \mlabel{rem:2.2} 
(a) In quantum mechanics the space 
$\bP(\cH^\infty)$ is interpreted as the state space of a physical system 
and the selfadjoint operator $-i\dd\pi(x)$ represents an observable. 
If $P$ is the spectral measure of this operator, then 
$-i\dd\pi(x) = \int_\R t\, dP(t)$ and 
\[ \Phi_\pi([v])(x) = \frac{1}{\la v,v\ra}
\int_\R t \, d\la P(t)v,v\ra \] 
is the expectation value of the observable $-i\dd\pi(x)$ in the state $[v]$.

(b) The groups arising as symmetry groups in quantum physics have natural 
projective representations (as symmetry groups of the quantum state space) 
satisfying a suitable positive energy condition (for the observable corresponding 
to  the Hamiltonian, resp.\ the energy). 
Accordingly, one often 
has to pass to central extensions to find unitary representations. 
\end{remark}

From \cite[Thm.~3.1]{Ne09} and its proof we obtain the following characterization 
of bounded representations: 
\begin{theorem}\mlabel{thm:2.3}
For a smooth unitary representation $(\pi,\cH)$, the following 
are equivalent: 
\begin{itemize}
\item[\rm(i)] $\pi$ is bounded. 
\item[\rm(ii)] $\pi \: G \to \U(\cH)$ is smooth with 
respect to the Banach--Lie group structure on $\U(\cH)$ defined by 
the norm topology. 
\item[\rm(iii)] $\cH^\infty = \cH$ and 
$\dd\pi \: \g \to B(\cH)$ is a continuous linear map. 
\item[\rm(iv)] $I_\pi$ is equicontinuous.
\end{itemize}
\end{theorem}

\section{Smoothing operators} 
\mlabel{sec:3}

A key tool to construct $C^*$-algebras for smooth unitary 
representations is the concept of a smoothing operator introduced 
in \cite{NSZ15}.

\begin{definition}
For a smooth unitary representation 
$(\pi, \cH)$ we call a bounded operator $A \in B(\cH)$ 
{\it smoothing} if $A\cH \subeq \cH^\infty$.\begin{footnote}{
In  \cite{KKW15} $A \in B(\cH)$ is called a  {\it Schwartz operator} if 
all operators $\dd\pi(D_1) A \dd\pi(D_2)$, $D_1, D_2 \in U(\g)$ are bounded. 
In view of Theorem~\ref{thm:smoothop}, for every 
Schwartz operator $A \in B(\cH)$ both $A$ and $A^*$ are smoothing. 
It is an interesting open problem whether the converse is also true. 
In \cite{KKW15} Schwartz operators are studied 
for the Schr\"odinger representation of the Heisenberg group on $\R^{2N}$. 
}\end{footnote}
\end{definition}

\begin{remark} (a) $\1 = \id_\cH$ is smoothing if and only if 
$\pi$ is bounded (cf.\ Theorem~\ref{thm:smoothop} below). 

(b)  If $G$ is finite dimensional and $f \in C^\infty_c(G)$, then 
$\pi(f) = \int_G f(g)\pi(g)\, dg$ is smoothing 
(\cite{Ga47}). 

Other types of smoothing operators can be obtained 
from a basis $x_1, \ldots, x_n$ of $\g$. 
{\it Nelson's Laplacian}
$\Delta := \oline{\sum_{j = 1}^n \dd\pi(x_j)^2}$ is selfadjoint with 
non-positive spectrum, so that the operators $(e^{t\Delta})_{t > 0}$ 
(the corresponding ``heat semigroup'') 
define a one-parameter group of selfadjoint contractions. 
It follows from \cite[Cor.~9.3]{Nel69} that these are smoothing operators. 

(c) If $\beta \: H \to G$ is a smooth homomorphism of Lie groups 
and $H$ is finite dimensional with 
$\cH^\infty = \cH^\infty(\pi\circ \beta)$, 
then, for every $f \in C^\infty_c(H)$, 
$\pi(f) = \int_H f(h)\pi(\beta(h))\, dh$ is smoothing by (b).
\end{remark} 

The following theorem is the main result on smoothing operators 
in \cite{NSZ15}. Its main power 
lies in the fact that the rather weak smoothing condition 
implies the smoothness of the multiplication maps in the operator norm. 

\begin{theorem} \mlabel{thm:smoothop} {\rm(Characterization Theorem for smoothing operators)} 
Let $(\pi, \cH)$ be a smooth unitary representation of a 
metrizable Lie group~$G$ with exponential 
function. Then $A \in B(\cH)$ 
is smoothing if and only if $G \to B(\cH), g \mapsto \pi(g) A$ 
is smooth with respect to the norm topology on $B(\cH)$. 
If, in addition, $G$ is Fr\'echet, then this is also equivalent to 
\begin{itemize}
\item[\rm(i)] $A\cH \subeq \cD(\oline{\dd\pi}(x_1) \cdots \oline{\dd\pi}(x_n))$ 
for $x_1, \ldots, x_n \in \g, n \in \N$. 
\item[\rm(ii)] All operators $A^* \oline{\dd\pi}(x_1) \cdots \oline{\dd\pi}(x_n)$, 
$x_1, \ldots, x_n \in \g, n \in \N$, 
are bounded. 
\end{itemize}
\end{theorem}

The following theorem is an important source of smoothing operators 
(\cite[Thm.~3.4]{NSZ15}): 
\begin{theorem}[Zellner's Smooth Vector Theorem]\mlabel{thm:3.3}
If $(\pi, \cH)$ is semibounded and $x_0 \in W_\pi$, 
then $\cH^\infty = \cH^\infty(\pi_{x_0})$ (cf.~\eqref{eq:stone}). 
In particular, the operators 
\[ e^{i\oline{\dd\pi}(x_0)} \quad \mbox{ and } \quad 
\int_\R f(t) \pi(\exp t x_0)\, dt = \hat f(i\oline{\dd\pi}(x_0)), 
\quad f \in \cS(\R), \]
are smoothing. 
\end{theorem}

The essentially selfadjoint operator  
$i\dd\pi(x_0)$ plays for a semibounded representation a similar role 
as Nelson's Laplacian $\Delta$ for a representation of a finite dimensional 
Lie group. This clearly demonstrates that, 
although one has very general 
tools that work for all representations of finite dimensional 
Lie groups (such as Nelson's hear semigroup and smoothing by convolution), 
specific classes of representations of infinite dimensional groups 
(such as semibounded ones) require specific but nevertheless equally powerful methods 
(such as Zellner's smoothing operators). 

\section{$C^*$-algebras} 
\mlabel{sec:4}

In the unitary representation theory of a finite dimensional 
Lie group $G$, a central tool is the convolution 
algebra $L^1(G)$, resp., its enveloping $C^*$-algebra $C^*(G)$, 
whose construction is based on the Haar measure, 
whose existence follows from the 
local compactness of $G$. Since the non-degenerate representations 
of $C^*(G)$ are in one-to-one correspondence with continuous unitary 
representations of $G$, the full power of the rich theory of $C^*$-algebras 
can be used to study unitary representations of $G$ (\cite{Dix77}). 

For infinite dimensional Lie groups 
there is no immediate analog of the convolution algebra 
$L^1(G)$, so that we cannot hope 
to find a $C^*$-algebra whose representations are in one-to-one 
correspondence with all unitary representations of $G$. However, 
in \cite{Gr05} H.~Grundling introduced the notion of a {\it host algebra}  
of a topological group $G$ and we shall see below that 
this concept provides 
natural $C^*$-algebras whose representations are in one-to-one 
correspondence with certain semibounded representations of~$G$.

\begin{definition} 
Let $\cA$ be a $C^*$-algebra, $M(\cA)$ its multiplier algebra,
\begin{footnote}{If $\cA$ is realized as a closed subalgebra of some 
$B(\cH)$, then 
\[ M(\cA) \cong \{ B \in B(\cH)\:(\forall A \in \cA)\ BA, AB \in \cA\}.\]}
\end{footnote}
and let $G$ be a topological group. 
We consider a homomorphism $\eta \: G \to \U(M(\cA))$ into the unitary 
group of the $C^*$-algebra $M(\cA)$. Then the pair $(\cA,\eta)$ is called 
a {\it host algebra for $G$} if,  
for each non-degenerate representation $\rho$ of ${\cal A}$ 
and its canonical extension $\tilde\rho$ to $M({\cal A})$, 
the unitary representation $\rho_G := \tilde\rho \circ \eta$ 
of $G$ is continuous and determines $\rho$ uniquely. 

In this sense, ${\cal A}$ is hosting a 
certain class of representations of $G$. 
We identify the set $\hat\cA$ of equivalence classes 
of irreducible representations of $\cA$ via 
$[\rho] \mapsto [\rho_G]$ with a subset $\hat G_\cA$ of the {\it unitary dual} $\hat G$ of $G$ 
(the set of equivalence classes of irreducible unitary representations of $G$). 
A host algebra $\cA$ is called 
{\it full} if it is hosting all continuous unitary representations of $G$.  
\end{definition}

If $G$ is locally compact, then $\cA = C^*(G)$, with 
$\eta$ specified by \break 
\[ \eta(g)f = \delta_g * f \quad \mbox{  for } \quad f \in L^1(G), g \in G,\]  
defines a full host algebra \cite{Dix77}. 
But to which extent do infinite dimensional 
Lie groups possess host algebras?  
If $G = (E,+)$ is an infinite dimensional locally convex space, 
then the set of equivalence classes of irreducible unitary representations 
can be identified with the dual space $E'$, and since this space carries
no natural locally compact topology, one cannot expect the existence 
of a full host algebra in general. 
Therefore one is looking for host algebras that 
accommodate certain classes of unitary representations. 
For bounded representations this is easy: 

\begin{example} \mlabel{ex:4.2} 
Let $(\pi, \cH)$ be a bounded unitary representation 
of the Lie group  $G$ and consider the 
$C^*$-subalgebra $\cA_\pi := C^*(\pi(G)) \subeq B(\cH)$ generated by 
$\pi(G)$. Then $\cA_\pi$ is unital, so that $\cA_\pi \cong M(\cA_\pi)$ and
$\eta := \pi \: G \to \U(\cA_\pi)$ is smooth with respect to the norm 
topology on the Banach--Lie group $\U(\cA_\pi)$. 
This implies that, for every representation $(\rho,\cK)$ 
of $\cA_\pi$, the representation $\rho_G := \rho \circ \pi$ is bounded 
and, since $\Spann \pi(G)$ is dense in $\cA_\pi$, it determines 
$\rho$ uniquely. Therefore $(\cA_\pi,\pi)$ is a host algebra 
of $G$ whose representations correspond to certain bounded representations 
of $G$. 

From this observation and Theorem~\ref{thm:2.3} one easily 
obtains for every continuous seminorm $p$ on $\g$ a 
a unital host algebra $(\cA_p, \eta_p)$ whose 
representations are precisely those smooth representations 
$(\pi,\cH)$ of $G$ satisfying $\|\dd\pi(x)\| \leq p(x)$ for $x \in \g$ 
(cf.\ \cite[\S III.2]{Ne00}), and this condition is equivalent to 
\[ I_\pi \subeq \{ \lambda \in \g' \: (\forall x \in \g)\, 
|\lambda(x)| \leq p(x)\}, \] 
which is an equicontinuous subset of~$\g'$. 
\end{example}

It turns out that Theorem~\ref{thm:3.3} is precisely what is needed 
to generalize the preceding construction to semibounded 
representations. Here one has to deal with non-unital 
algebras and smoothing operators of the form $e^{i\oline{\dd\pi}(x_0)}$ 
that lead for every semibounded representation 
$(\pi,\cH)$ and $x_0 \in W_\pi$ to the host algebra 
$\cA := C^*\big(\pi(G) e^{i\oline{\dd\pi}(x_0)} \pi(G)\big)$. 
Putting everything together, we obtain: 

\begin{theorem} {\rm(\cite[Cor.~4.9]{NSZ15})}  \mlabel{thm:5.2} 
Let $C \subeq \g'$ be a weak-$*$-closed $\Ad^*(G)$-invariant subset
which is semi-equicontinuous in the sense that its support function 
$s_C(x) := \sup \la C, x \ra$
is bounded in a neighborhood of some $x_0 \in \g$. Then there exists for every 
semibounded representation $(\pi, \cH)$ a 
host algebra $(\cA_C,\eta_C)$ of $G$ whose representations correspond to those 
semibounded unitary representations $(\pi, \cH)$ of $G$ for which 
$s_\pi \leq s_C$, resp., $I_\pi \subeq -C$. 
\end{theorem}

\begin{remark} (\cite[Prop.~6.13]{Ne08})\mlabel{rem:4.4}
(a) If $\pi$ is (semi)bounded, then $I_\pi$ is a (locally) compact subset 
of $\g'$, endowed with the weak-$*$-topology and 
for every $x_0 \in W_\pi$ the map 
$I_\pi \to \R, \alpha \mapsto \alpha(x_0)$ is proper. 

(b) If $G = (E,+)$ is the additive group of a locally convex space, then 
we identify the topological dual space $E'$ with the character group 
$\hat G$ by $\chi_\alpha(v) := e^{i\alpha(v)}$. Then there exists for 
every semibounded representation $(\pi, \cH)$ a spectral 
measure $P$ on the locally compact space $I_\pi \subeq E'$ with 
$\pi(v) = \int_{E'} e^{i\alpha(v)}\, dP(\alpha)$ and the $C^*$-algebra 
from Theorem~\ref{thm:5.2} is isomorphic to 
$C_0(\supp(P))$, where $\supp(P) \subeq I_\pi$ is the support of~$P$ 
(cf.\ \cite[Thm.~4.1]{Ne09}). 

(c) If $C_1 \subeq C_2$ are $\Ad^*(G)$-invariant weak-$*$-closed convex 
equicontinuous subsets, then the construction of the host algebras 
$\cA_{C_1}$ and $\cA_{C_2}$ provides a morphism 
$\cA_{C_2} \onto \cA_{C_1}$. Since the collection of all 
$\Ad^*(G)$-invariant weak-$*$-closed convex equicontinuous subsets 
is directed with respect to inclusion, we thus obtain a projective 
system of $C^*$-algebras $\cA_C$ whose representation theory 
completely describes the semibounded representations of~$G$. 
Similar results for finite dimensional groups obtained by 
holomorphic extensions to complex semigroups 
can already be found in \cite[\S XI.6]{Ne00}. 
\end{remark}

\begin{remark} In general, continuous unitary representations of 
locally convex spaces can not be represented in terms of spectral 
measures on $E'$. This is closely related to the problem of writing 
the continuous positive definite functions 
$\pi^{v,v}(x) := \la \pi(x)v,v\ra$ as a Fourier transform  
$$ \hat\mu(x) = \int_{E'} e^{i\alpha(x)}\, d\mu(\alpha) $$
of some finite measure $\mu$ on $E'$. If $E$ is nuclear, 
then the Bochner--Minlos Theorem \cite{GV64} ensures the existence 
of such measures and hence also of spectral measures representing 
unitary representations. However, if $E$ is an infinite dimensional 
Hilbert space $E$, then the continuous positive definite 
function $\phi(v) := e^{-\|v\|^2/2}$ is not the Fourier transform 
of a positive measure on $E'$ (\cite[Ex.~17.1]{Ya85}). 
Note that $E$ is not nuclear, so that 
the Bochner--Minlos Theorem does not apply. 
Therefore it is quite remarkable that 
nuclearity assumptions are not needed to deal with semibounded 
representations because the domains of the spectral measures are locally compact. 
It is easy to show that, if $E$ is separable and metrizable, then a spectral measure 
for a unitary representation 
$(\pi,\cH)$ exists if and only if $\pi$ is a direct sum of bounded representations. 
\end{remark}

The main motivation to find host algebras 
is that they provide very well developed tools to 
decompose the representations into irreducible ones: 
If $(\cA,\eta)$ is a separable host algebra of $G$ which is of type I,\begin{footnote}
  {This means that, for every (non-zero) irreducible representations 
$(\rho,\cK)$ of $\cA$, the image $\rho(\cA)$ contains 
a non-zero compact operator (\cite{Sa66}), and this in turn implies 
that $K(\cK) \subeq \rho(\cA)$.}
\end{footnote}
then the following 
abstract disintegration theorem applies immediately to all 
$\cA$-representations $\rho_G$ of $G$ 
and thus reduces the classification problems for semibounded 
representations to the classification of the irreducible representations 
and the spectral multiplicity theory of measures on~$\hat G_\cA$:

\begin{theorem}[Abstract Disintegration Theorem] 
{\rm(\cite[Th.\ 8.6.6]{Dix77})} 
Let $\cA$ be a separable type I $C^*$-algebra. 
\begin{itemize}
\item[\rm(i)] Every Borel measure $\mu$ on $\hat\cA$ defines a 
multiplicity free direct integral representation 
$\pi_\mu := \int_{\hat \cA}^\oplus \pi\, d\mu(\pi).$ 
Two such representations $\pi_\mu$ and $\pi_\nu$ are equivalent if and only 
the measure classes of $\mu$ and $\nu$ coincide. 
\item[\rm(iii)] For every separable representation $(\pi, {\cal H})$,
there exist mutually disjoint measures $(\mu_n)_{n \in
\N\cup\{\infty\}}$ such that 
\[ \pi \cong \pi_{\mu_1} \oplus 2 \cdot\pi_{\mu_2} \oplus \cdots \oplus
\aleph_0 \cdot\pi_{\mu_\infty}.
\begin{footnote}
{For a representation $(\rho, {\cal K})$ and a cardinal $n$, we write 
$n\cdot \rho$ for the representation $\rho \otimes \1$ on 
${\cal K} \hat\otimes \ell^2(\bn)$, where $\bn$ is a set of cardinality $n$.}  
\end{footnote}
\] 
The measure classes of $(\mu_n)_{n \in \N \cup
\{\infty\}}$ are uniquely determined by $\pi$. 
\end{itemize}
\end{theorem}

\section{Representations of unitary groups} 
\mlabel{sec:5}

Let $\cH$ be an infinite dimensional complex Hilbert space. 
We write $\U(\cH)_n$ for the Banach--Lie group obtained by endowing the 
unitary group $\U(\cH)$ with the norm topology and 
$\U(\cH)_s$ for the topological group structure obtained from the 
strong operator topology, i.e., the topology of pointwise convergence.

The representation theory of infinite dimensional unitary groups began with 
I.~E.~Segal's paper \cite{Se57}, where he studies so-called 
{\it physical representations} of the full group $\U({\cal H})_s$.
These are characterized by the
condition that their differential maps finite rank hermitian
projections to positive operators. Segal shows that 
physical representations are direct sums of irreducible ones, 
which are precisely those 
occurring in the decomposition of finite tensor products 
$\cH^{\otimes N}$, $N \in\N_0$. 
This tensor product decomposes as in classical 
Schur--Weyl theory: 
\begin{equation}
  \label{eq:schur-weyl}
\cH^{\otimes N} \cong \bigoplus_{\lambda \in \Part(N)} \bS_\lambda(\cH) \otimes \cM_\lambda,
\end{equation}
where $\Part(N)$ is the set of all partitions 
$\lambda = (\lambda_1, \ldots, \lambda_n)$ of $N$, and 
$\bS_\lambda(\cH)$ is an irreducible unitary representation of $\U(\cH)$ 
(called a {\it Schur representation}), and 
$\cM_\lambda$ is the corresponding irreducible representation of the symmetric group $S_N$, 
hence in particular finite dimensional.\begin{footnote}
{We refer to \cite{BN12} for an extension of Schur--Weyl theory to irreducible 
representations of $C^*$-algebras.} 
\end{footnote}
In particular, 
$\cH^{\otimes N}$ is a finite sum of irreducible representations of~$\U(\cH)$. 

The representation theory of the Banach--Lie group  
\[ \U_\infty({\cal H})  = \U(\cH) \cap (\1 + K(\cH)),\]
 where ${\cal H}$ is a separable 
complex Hilbert space and $K(\cH)$ the space of compact 
operators, was developed by Kirillov and Olshanski 
in \cite{Ki73} and \cite[Thm.~1.11]{Ol78}. They show that 
all continuous representations of $\U_\infty(\cH)$
are direct sums of irreducible 
representations and that, for $\K = \C$, all irreducible 
representations are of the form $\bS_\lambda(\cH) \otimes \bS_\mu(\oline\cH)$, where 
$\oline\cH$ is the space $\cH$, endowed with the opposite complex structure. 
They also obtained generalizations 
 for the corresponding groups over real and quaternionic 
Hilbert spaces. It follows in particular that all irreducible
representations $(\pi, \cH_\pi)$ of the Banach--Lie group $\U_\infty(\cH)$ are 
bounded. The classification of the bounded unitary representations 
of the Banach--Lie group $\U_p(\cH) := \U(\cH) \cap (\1 + B_p(\cH))$ 
(where $B_p(\cH)$ is the $p$th Schatten ideal) 
does not depend on $p$ as long as $1 < p \leq \infty$, 
but, for $p = 1$, factor representations of type II and III exist 
(see \cite{Boy80} for $p = 2$, and \cite{Ne98} for the general case). 
Dropping the boundedness assumptions leads to non-type I 
representations of $\U_p(\cH)$, $p < \infty$ 
(cf.\ \cite[Thm.~5.5]{Boy80}). 

These results clearly show that the group $\U_\infty(\cH)$ is singled out among all its relatives $\U_p(\cH)$  
by the fact that its unitary representation theory is well-behaved. 
If $\cH$ is separable, then $\U_\infty(\cH)$ is separable, 
so that its cyclic representations are separable as well. 
Hence there is no need to discuss inseparable representations for this group. 
This is different for the Banach--Lie group 
$\U(\cH)_n$ which has many inseparable bounded irreducible unitary representations coming 
from irreducible representations of the Calkin algebra $B(\cH)/K(\cH)$. 
The following theorem was an amazing achievement 
of D.~Pickrell \cite{Pi88}, showing 
that restricting attention to 
representations on separable spaces tames the representation theory 
of $\U(\cH)_n$ in the sense that all its separable representations 
are actually continuous with respect to the strong operator topology, i.e., continuous 
representations of $\U(\cH)_s$. 

\begin{theorem}[Pickrell's Continuity Theorem] \mlabel{thm:pickrell} 
Every separable unitary representation
\begin{footnote}{Recall that we included continuity in the definition 
of a unitary representation.} \end{footnote}
of $\U(\cH)_{n}$ is also 
continuous for the strong operator topology,  hence a representation of $\U(\cH)_{s}$. 
\end{theorem}

Building on the fact that the identity component 
$\U_\infty(\cH)_0$ is dense in $\U(\cH)_s$ 
and extending the Kirillov--Olshanski classification to non-separable Hilbert spaces, 
we show in \cite{Ne14b} that the representations of 
$\U_\infty(\cH)_0$ and $\U(\cH)_s$ coincide, more precisely: 

\begin{theorem}
Let $\cH$ be an infinite dimensional real, complex or quaternionic 
Hilbert space. Then every continuous unitary representation of $\U_\infty(\cH)_0$ 
extends uniquely to a continuous representation of $\U(\cH)_s$. 
All unitary representations of these groups are direct sums of irreducible ones 
which are of the form \break $\bS_\lambda(\cH) \otimes \bS_\mu(\oline\cH)$, hence in particular 
bounded representations of $\U(\cH)_n$. 
\end{theorem}

In view of the preceding results, the separable representation theory 
of the Lie group $\U(\cH)_n$ very much resembles the representation 
theory of a compact group, and so does the representation theory of $\U_\infty(\cH)$. 

\section{Bounded representations of gauge groups} 
\mlabel{sec:6}

Let $M$ be a smooth $\sigma$-compact manifold, $K$ a compact connected Lie group and 
$q \: P \to M$ a $K$-principal bundle. 
We consider the Lie group $G := \Gau_c(P)$ 
of compactly supported gauge transformations and observe that it 
can be realized as 
\[ \Gau_c(P) \cong C^\infty_G(P,K) := \{ f \in C^\infty_c(P,K) \: 
(\forall p \in P)(\forall k \in K)\ f(p.k) = k^{-1}f(p)k\}.\] 
Its Lie algebra $\g$ is 
\[ \gau_c(P) \cong C^\infty_G(P,\fk) := \{ \xi \in C^\infty_c(P,\fk) \: 
(\forall p \in P)(\forall k \in K)\ \xi(p.k) = \Ad(k)^{-1}\xi(p)\}.\] 

Bounded representations of $G$ are easy to construct by evaluations. 
To see this,  let $x \in P$, $(\rho,V_\rho)$ be an irreducible representation of $K$ 
(which is automatically finite dimensional). Then 
\[ \pi_{x,\rho}(f) := \rho(f(x)) \quad \mbox{ for } \quad f\in G, x\in P \] 
defines a finite dimensional irreducible representation of $G$. 
Clearly $\pi_{x,\rho} \cong \pi_{y,\rho}$ if $q(x) = q(y)$. 
For finite subsets 
$\bx = \{x_1, \ldots,x_N\} \subeq P$ for which the points 
$q(x_i) \in M$ are pairwise different, the representation 
\[ \pi_{\bx, \rho}:= \bigotimes_{x \in \bx} \pi_{x, \rho_x} \] 
of $G$ is also bounded and irreducible  
\cite{JN14}. 

We can even go further: Let $\bx \subeq P$ be a subset mapped injectively by 
$q$ to a locally finite subset of $M$. Then we assign to every collection 
$\rho = (\rho_x, V_x)_{x \in \bx}$ of irreducible representations of $K$ 
the UHF-$C^*$-algebra 
\[ \cA_{\bx,\rho} := \bigotimes_{x \in \bx} B(V_x). 
\quad \mbox{ Then }\quad 
 \eta_{\bx, \rho}:= \bigotimes_{x \in \bx} \pi_{x, \rho_x} \: 
C^\infty_c(M,K) \to \U(\cA_{\bx,\rho}) \] 
defines a smooth  homomorphism into the unitary group of the 
$C^*$-algebra $\cA_{\bx,\rho}$ and it is not hard to see that 
the image of $\eta_{\bx,\rho}$ generates $\cA_{\bx,\rho}$, so that 
we obtain a host algebra $(\cA_{\bx, \rho}, \eta_{\bx,\rho})$ of $G$ 
(cf.~Example~\ref{ex:4.2}). The following theorem reduces the classification 
of the bounded representations of $G$ completely to questions on 
$C^*$-algebras (cf.\ \cite{JN14}). 

\begin{theorem} 
Every bounded irreducible representation $\pi$ of the identity component 
$\Gau_c(P)_0$ is of the form 
$\beta \circ \eta_{\bx,\rho}$ for some irreducible 
representation $\beta$ of $\cA_{\bx,\rho}$. 
If $M$ is compact, then $\bx$ is {finite} and $\pi \cong \pi_{\bx,\rho}$. 
\end{theorem}

\begin{remark} (a) If the bundle $P$ is trivial, then 
$G \cong C^\infty_c(M,K)$ is a mapping group. 
It is contained in the Banach--Lie group 
$C_0(M,K)_0$ of maps vanishing at infinity. 
In \cite{JN14} it is also shown that every irreducible 
bounded representation of $G$ that extends to a bounded 
representation of this group is of the form 
$\pi_{\bx,\rho}$ for a finite subset $\bx$. 

(b)  A description of the irreducible representations of 
the UHF-$C^*$-algebras 
$\cA = \bigotimes_{n \in \N} M_{d_n}(\C)$ can be obtained from the work 
of Glimm \cite{Gli60} and Powers \cite{Po67}. 
One of their main results is that 
all irreducible representations of $\cA$ are twists of 
infinite tensor products of irreducible representations 
$\bigotimes_x (V_x, v_x)$ ($v_x \in V_x$ a unit vector) 
by an automorphism of $\cA$.  Since $\cA$ has a large automorphism group, 
this implies in particular that 
not all bounded irreducible representations $G$ are tensor products 
of irreducible representations of the factors. 
\end{remark}

\section{From bounded to semibounded representations by 
holomorphic induction} 
\mlabel{sec:7}

We now describe a complex geometric technique to classify semibounded representations. 
For the results in this section we refer to \cite[App.~C]{Ne14} for 
the Fr\'echet case and to \cite{Ne13} for a more 
detailed exposition of these methods for Banach--Lie groups. 

Let $G$ be a connected Fr\'echet--BCH--Lie group
with Lie algebra $\g$.\begin{footnote}{A Lie group $G$ is called {\it locally exponential}  
if $\exp \: \g \to G$ maps some $0$-neighborhood in $\g$ diffeomorphically 
to an open $\1$-neighborhood in $G$. If, in addition, $G$ is an 
{\it analytic Lie group} and the local diffeomorphism is bianalytic, 
then $G$ is called a {\it BCH--Lie group} because this implies 
that the Baker--Campbell--Hausdorff (BCH) series defines a local 
analytic group structure on some $0$-neighborhood in $\g$ 
(cf.~\cite[\S IV.1]{Ne06}). 
}\end{footnote}
We further assume that 
there exists a complex BCH--Lie group $G_\C$ with Lie algebra $\g_\C$ and a natural 
map $\eta \: G \to G_\C$ for which $\L(\eta)$ is the inclusion $\g \into \g_\C$. 
Let $H \subeq G$ be a Lie subgroup with Lie algebra $\fh$ 
for which $M := G/H$ carries the structure of a smooth manifold 
with a smooth $G$-action. 
We also assume the existence of closed $\Ad(H)$-invariant subalgebras 
$\fp^\pm \subeq \g_\C$ with $\oline{\fp^\pm} = \fp^\mp$ for $\oline{x+iy} := x- i y$, 
$x,y\in \g$, and for which we have a topological direct 
sum decomposition 
\begin{equation}
  \label{eq:splitcond}
 \g_\C = \fp^+ \oplus \fh_\C \oplus \fp^-.\tag{SC}
\end{equation}
We put 
\[ \fq := \fp^+ \rtimes \fh_\C \quad \mbox{ and } \quad \fp := \g \cap (\fp^+ \oplus \fp^-),\]
so that $\g = \fh \oplus \fp$ is a topological direct sum. We 
assume that there exist open symmetric convex $0$-neighborhoods 
$U_{\g_\C} \subeq \g_\C$, 
$U_\fp \subeq \fp \cap U_{\g_\C}, U_\fh \subeq \fh \cap U_{\g_\C}, U_{\fp^\pm} \subeq\fp^\pm  \cap U_{\g_\C}$ and  
$U_\fq \subeq \fq \cap U_{\g_\C}$ such that the 
BCH-product $x * y = x + y + \frac{1}{2} [x,y] + \ldots$ is defined and holomorphic on 
$U_{\g_\C} \times U_{\g_\C}$, and the following maps are analytic diffeomorphisms onto an open subset: 
\begin{description}
\item[\rm(A1)] $U_{\fp} \times U_\fh \to \g, (x,y) \mapsto x * y$. 
\item[\rm(A2)] $U_{\fp} \times U_\fq \to \g_\C, (x,y) \mapsto x * y$. 
\item[\rm(A3)] $U_{\fp^-} \times U_\fq \to \g_\C, (x,y) \mapsto x * y$. 
\end{description}

Then (A1) implies the existence of a smooth manifold structure on $M = G/H$ on 
which $G$ acts analytically. 
Condition (A2) implies the existence of a complex manifold structure on $M$ which is $G$-invariant 
and for which $T_{\1 H}(M) \cong \g_\C/\fq$. Finally, 
(A3) makes the proof of \cite[Thm.~2.6]{Ne13} work, so that we can associate to every 
bounded unitary representation $(\rho,V)$ of $H$ 
a holomorphic Hilbert bundle $\bV := G \times_H V$ over the complex 
$G$-manifold $M$.  

\begin{definition}
 \mlabel{def:d.1} 
We write $\Gamma(\bV)$ for the space of holomorphic sections 
of the holomorphic Hilbert bundle $\bV \to M = G/H$ on which the group $G$ acts by 
holomorphic bundle automorphisms. 
A unitary representation $(\pi, \cH)$ of $G$ is said to be 
{\it holomorphically induced from $(\rho,V)$} 
if there exists a $G$-equivariant linear injection 
$\Psi \: \cH \to \Gamma(\bV)$ such that the 
adjoint of the evaluation map $\ev_{\1 H} \: \cH \to V = \bV_{\1 H}$ 
defines an isometric embedding $\ev_{\1 H}^* \: V \into \cH$. 
If a unitary representation $(\pi, \cH)$ holomorphically induced 
from $(\rho,V)$ exists, then it is uniquely determined 
(\cite[Def.~3.10]{Ne13}) and we call $(\rho,V)$ {\it (holomorphically)  
inducible}. 

This concept of inducibility involves a choice of sign. 
Replacing $\fp^+$ by $\fp^-$ changes the complex structure on $G/H$ 
and thus leads to a different class of holomorphically inducible 
representations of of $H$. 
\end{definition}

\begin{theorem} \mlabel{thm:c.3}
Suppose that $(\pi, \cH)$ is a unitary representation of $G$ and 
$V \subeq \cH$ is an $H$-invariant closed subspace such that 
\begin{description}
\item[\rm(HI1)] the representation $(\rho,V)$ of $H$ on $V$ is bounded,  
\item[\rm(HI2)] $V \cap (\cH^\infty)^{\fp^-}$ is dense in $V$,  and 
\item[\rm(HI3)] $\pi(G)V$ spans a dense subspace of $\cH$. 
\end{description}
Then $(\pi, \cH)$ is holomorphically induced from $(\rho,V)$ 
and $\pi(G)' \to \rho(H)', A \mapsto A\res_V$ is an isomorphism of the commutants. 
\end{theorem}

The preceding theorem implies in particular that $(\rho,V)$ is irreducible 
if and only if $(\pi, \cH)$ is. All the concrete classification 
results for semibounded irreducible representations rest on the 
fact that they can be obtained by bounded representations of suitably 
chosen subgroups $H$ for which a classification of the irreducible 
bounded representations is available. 

\begin{example} It is instructive to see how the general method of holomorphic 
induction matches the classification of irreducible unitary representations
 of a compact connected Lie group $G$. 
In this case we choose a maximal torus $H \subeq G$ and obtain a triangular 
decomposition as in \eqref{eq:splitcond}, where 
$\fb = \fp^- \rtimes \fh_\C$ 
is a Borel subalgebra of the complex reductive Lie algebra $\g_\C$. 
Then the bounded irreducible representations of 
$H$ are one-dimensional, hence  given by characters $\chi \: H \to \T$.
Such a character is holomorphically inducible if and only if 
the weight $\dd\chi \: \fh_\C\to \C$ is $\fb$-dominant. We thus 
obtain the well-known classification of the finite dimensional irreducible 
representations of $G$ by $\fb$-dominant weights on $\fh_\C$. 
\end{example} 

\section{Hermitian groups} 
\mlabel{sec:8}

We now explain the main points of the classification 
of irreducible semibounded representations of hermitian 
Lie groups carried out in \cite{Ne12}. 

\begin{definition}
(a) A {\it hermitian Lie group} is a triple $(G,\theta, \bd)$, where 
$G$ is a connected Lie group, $\theta$ an involutive 
automorphism of $G$ with the corresponding eigenspace 
decomposition $\g = \fk \oplus \fp$, 
$\bd \in \z(\fk)$ (the center of $\fk$) 
an element for which $\ad \bd\res_\fp$ is a complex 
structure, and $\fp$ carries an $e^{\ad \fk}$-invariant 
Hilbert space structure. 
We then write $K := (G^\theta)_0$ for the identity component of 
the group of $\theta$-fixed points in $G$ and observe that our 
assumptions imply that $G/K$ is a hermitian symmetric space 
(modeled on a complex Hilbert space).   

(b) We call $(G,\theta,\bd)$ {\it irreducible} 
if the unitary $K$-representation on $\fp$ is irreducible. 
We say that $\g$ is {\it full} if 
$\ad \fk \subeq \fu(\fp^+)$ is the full derivation algebra 
with respect to the Jordan product $[x,y,z] := [[x,\oline y],z]$ 
on the $i$-eigenspace $\fp^+ \subeq \g_\C$. 
\end{definition}

If $\g$ is full, then the Lie algebra $\fk/\fz(\fk)$ contains 
no non-trivial open invariant cones (\cite[Lemma~5.10]{Ne12}). 
If, in addition, $(G,\theta,\bd)$ is irreducible 
and $(\pi, \cH)$ an irreducible semibounded representation, 
then either $\bd \in W_\pi$, i.e., $\pi$ is a {\it positive energy representation}, 
or the dual representation $\pi^*$ has this property 
(\cite[Thm.~6.2]{Ne12}). In the former case, it is holomorphically 
induced from a bounded representation $(\rho,V)$ of $K$ 
(\cite[Thm.~6.4]{Ne12}). Therefore the classification of irreducible 
semibounded representations completely reduces the determination of the 
irreducible bounded representations of~$K$ which are holomorphically 
inducible. 

\begin{example}  \mlabel{ex:8.2} If $G = \fp \rtimes_\alpha K$ is 
a Cartan motion group of the complex Hilbert space $\fp$, i.e., 
$[\fp,\fp] = \{0\}$, then all semibounded representations 
$(\pi, \cH)$ of $G$ are trivial on $\fp$ by 
\cite[Thm.~7.1]{Ne12}. 

However, the Lie algebra has a unique central extension 
$\hat \g = \R \oplus_\omega \g$, 
given by the cocycle $\omega(x,y) := \Im \la x_\fp,y_\fp\ra$ for 
$x = x_\fk + x_\fp$. 
If $\g$ is full, then $\hat K \cong \R \times \U(\fp)$ and 
$(\rho,V)$ is holomorphically inducible if and only if 
$-i\dd\rho(1,0) \geq 0$ 
(\cite[Thm.~7.2]{Ne12}). 
\end{example}

\begin{example} \mlabel{ex:8.3} 
(a) For the case where $(G,\theta,\bd)$ is irreducible and full and 
$\cD := G/K$ is an infinite dimensional
 symmetric Hilbert domain (a generalization of a bounded symmetric domain in $\C^n$), 
then, up to representations vanishing on $\fp$, the irreducible 
semibounded representations of $G$ are determined in \cite[Thm.~8.3]{Ne12}. 
Here it is natural to assume that $G$ is the universal central extension 
of the connected automorphism group of $\cD$. Then $K$ is a product of at most 
three factors isomorphic to $\R$, an infinite dimensional group 
$\U(\cH)_n$, or to a simply connected covering group $\tilde\U_n(\C)$. 
Therefore the separable irreducible bounded representations of $K$ are 
well-known (cf.\ Section~\ref{sec:5}), 
so that the main point is to obtain the inequalities 
characterizing holomorphic inducibility. It is remarkable that, in all cases, the central 
extension is needed to have semibounded representations 
that are non-trivial on $\fp$. 

(b) A concrete example of a full hermitian group is the universal central 
extension $G = \hat \Sp_{\rm res}(\cH)$ 
of the restricted symplectic group of a complex Hilbert space 
with respect to the symplectic form $\omega(v,w) := \Im \la v, w \ra$: 
\[ \Sp_{\rm res}(\cH) := \{ g \in \Sp(\cH,\omega) \: \|gg^\top - \1\|_2 < \infty \}. \] 
Here $K \cong \R \times \U(\cH)$. An important example of a semibounded 
representation is the metaplectic representation 
of $\hat\Sp_{\rm res}(\cH)$ on the bosonic Fock space 
$S(\cH) = \hat\oplus_{n \in\N_0} S^n(\cH)$ (\cite{Se78}). 

(c) Another example of a full irreducible hermitian Lie group is 
the conformal group $G = \OO(\R^2,\cH)$ of an infinite dimensional Minkowski space. 
Neither $G$ nor any of its covering groups have non-trivial semibounbed representations 
(\cite[Thm.~8.5]{Ne12}). 
\end{example}

\begin{example}  \mlabel{ex:8.4} 
(a) For the hermitian groups $G$ for which the Cartan dual Lie algebra 
 $\g^c = \fk \oplus i\fp$ is of the type considered in Example~\ref{ex:8.3}, 
i.e., $G/K$ is dual to a symmetric Hilbert domain, the classification of the 
semibounded representations is quite easy to describe because in this case 
a bounded representation $(\rho,V)$ of $K$ is holomorphically inducible 
if and only if $\rho$ is {\it anti-dominant} in the  sense that 
\[ \dd\rho([z^*,z]) \geq 0 \quad \mbox{ for }\quad z \in \fp^+ \subeq \g_\C\]  
(\cite[Thm.~9.1]{Ne12}). 

(b) A concrete example is the restricted unitary group of a 
Hilbert space $\cH = \cH_+ \oplus \cH_-$: 
\[ G = \U_{\rm res}(\cH,\cH_+) := \bigg\{\pmat{a & b \\ c & d}\in \U(\cH) \: \|b\|_2, \|c\|_2 < \infty\bigg\}, \]
for which $K \cong \U(\cH_+) \times \U(\cH_-)$ and 
$G/K$ is the {\it restricted Grassmannian} (cf.\ \cite{PS86}, \cite{Mick89}). 
\end{example}

\section{Loop groups (Affine Kac--Moody groups)} 
\mlabel{sec:9}

A {\it Hilbert--Lie algebra} is a real Lie algebra $\fk$ carrying the structure 
of a real Hilbert space, such that the scalar product $(\cdot,\cdot)$ is 
invariant under the adjoint action, i.e., 
\[ ([x,y],z) = (x,[y,z]) \quad \mbox{ for } \quad x,y,z\in \fk.\]
In finite dimensions, these are precisely the compact Lie algebras. In 
infinite dimensions, $\fk$ is a direct sum of an abelian ideal and simple 
Hilbert--Lie algebras which are isomorphic to 
$\fu_2(\cH)$ (the skew-hermitian Hilbert--Schmidt operators) on 
an infinite dimensional real, complex or quaternionic Hilbert space~$\cH$ 
(\cite{Sch60}). 

Let $K$ be a simply connected Lie group for which $\fk$ is a simple 
Hilbert--Lie algebra and $\phi \in \Aut(K)$ an automorphism with 
$\phi^N = \id_K$. Then the {\it twisted loop group} 
\[  \cL_\phi(K) 
:= \Big\{ f \in C^\infty(\R,K) \: (\forall t \in \R)\ 
f\Big(t + \frac{2\pi}{N}\Big) = \phi^{-1}(f(t))\Big\} \] 
is a Fr\'echet--Lie group with Lie algebra 
\[ \cL_\phi(\fk) 
:= \Big\{ \xi \in C^\infty(\R,\fk) \: (\forall t \in \R)\ 
\xi\Big(t + \frac{2\pi}{N}\Big) = \L(\phi)^{-1}(\xi(t))\Big\}. \] 
For $\phi = \id_K$, we obtain the loop group $\cL(K) = C^\infty(\bS^1,K)$ 
(see~\cite{PS86} for finite dimensional $K$ and \cite{Ne14} for the infinite 
dimensional case). 

The subgroup $\cL_\phi(K) \subeq C^\infty(\R,K)$ is 
translation invariant, so that we obtain by 
$\alpha_s(f)(t) := f(t + s)$ a smooth action 
of the circle group $\T \cong \R/2\pi \Z$ on $\cL_\phi(K)$. 
The cocycle 
$$ \omega(\xi, \eta) := \int_0^{2\pi} \la \xi'(t), \eta(t)\ra\, dt $$
defines a central extension 
$$\tilde\cL_\phi(\fk) := \R \oplus_\omega \cL_\phi(\fk), $$
which leads by $D(z,\xi) := (0,\xi')$ to the double extension 
$$\g := \hat\cL_\phi(\fk) := (\R \oplus_\omega \cL_\phi(\fk)) \rtimes_D \R. $$
We put $\bd := (0,0,-1) \in \g$.

To formulate the main results of \cite{Ne14}, 
we first recall that one has a complete classification of the twisted 
loop groups for infinite dimensional~$K$. There are four classes 
of loop algebras: the untwisted loop algebras  
$\cL(\fu_2(\cH))$, where $\cH$ is an infinite dimensional 
real, complex or quaternionic Hilbert space, and a twisted type 
$\cL_\phi(\fu_2(\cH))$, where $\cH$ is a complex Hilbert space 
and $\phi(x) = \sigma x \sigma$ holds for an antilinear isometric involution
$\sigma \: \cH \to \cH$ (this corresponds to complex conjugation of the corresponding 
matrices). We call $\phi$ of {\it standard type}\begin{footnote}
{This terminonology stems from the fact that the standard type automorphisms 
naturally realize the $7$ types of locally affine root systems 
$A_J^{(1)}, B_J^{(1)}, C_J^{(1)}, D_J^{(1)}, B_J^{(2)}, C_J^{(2)}$ and $BC_J^{(2)}$.}  
\end{footnote}
if either $\phi = \id_K$, or 
\begin{itemize}
\item $\K = \R$ and $\phi(g) := r g r^{-1}$, where $r$ is the orthogonal reflection in 
a hyperplane,
\item $\K = \C$, $\cH = \cH_0 \oplus  \cH_0$, 
$\sigma_0$ is a conjugation on $\cH_0$, 
$\sigma(x,y) := (\sigma_0(x), \sigma_0(y))$ on $\cH$, and 
\[ \phi(g) = S \sigma g \sigma  S^{-1} \quad \mbox{ for } \quad 
S = \pmat{ 0 & \1 \\ -\1 & 0}, \quad g \in K = \U_2(\cH).\]
\item $\K = \C$, $\cH = \cH_0 \oplus \C \oplus  \cH_0$, 
$\sigma(x,y,z) := (\sigma_0(x),\oline y, \sigma_0(z))$, and 
\[ \phi(g) = S \sigma g \sigma S^{-1} \quad \mbox{ for } \quad 
S = \pmat{ 0 & 0 & \1 \\ 0 & 1 & 0 \\ \1 & 0 & 0}, \quad g \in K = \U_2(\cH).\] 
\end{itemize}

In \cite{Ne14}, the classification 
of the semibounded unitary representations of the 
$1$-connected Lie group $G := \hat\cL_\phi(K)$ corresponding to the 
double extensions $\g = \hat\cL_\phi(\fk)$ is obtained if $\phi$ is of standard type. 
For a reduction of the general case to this one, we refer to 
\cite{MN15b}. The first major step is to show that, for an irreducible semibounded 
representation $(\pi, \cH)$, the operator 
$-i\dd\pi(\bd)$ is either bounded from below (positive energy representations) 
or from above (negative energy representations). 
Up to passing to the dual representation, we may therefore 
assume that we are in the first case. 
Then the minimal spectral value of $-i\dd\pi(\bd)$ turns out to be an eigenvalue 
and the centralizer $H := Z_G(\bd)$ of $\bd$ 
in $G$ acts on the corresponding eigenspace, 
which leads to a bounded irreducible representation $(\rho,V)$ of~$H$ 
from which $(\pi, \cH)$ is obtained by holomorphic induction. 
Since an explicit 
classification of the bounded irreducible representations 
of the groups $Z_G(\bd)_0$ is available from 
\cite{Ne98, Ne12} 
in terms of $\cW$-orbits of extremal weights, it remains to characterize 
those weights $\lambda$ for which the corresponding representation 
$(\rho_\lambda, V_\lambda)$ is holomorphically inducible. 
This is achieved in \cite[Thm.~5.10]{Ne14}. 
It is equivalent to $\lambda$ being $d$-minimal, and the final step 
consists in showing that the irreducible $G$-representation 
$(\pi_\lambda, \cH_\lambda)$ corresponding to a $d$-minimal 
weight is actually semibounded (\cite[Thm.~6.1]{Ne14}). 

For untwisted loop groups $\hat\cL(K)$ and compact groups $K$, 
the corresponding class of representations is well-known from 
the context of affine Kac--Moody algebras 
(cf.~\cite{PS86}). In this context one thus obtains 
the class of positive energy representations, 
but for infinite dimensional $K$ the positive energy 
condition is too weak to make holomorphic induction work 
(cf.\ \cite{Ne15}). 

\section{The diffeomorphism group of the circle} 
\mlabel{sec:10}

For the group $G := \Diff(\bS^1)_+$ of orientation preserving diffeomorphisms 
of the circle, the Lie algebra is the space 
\[ \g = \cV(\bS^1) = C^\infty(\bS^1) \frac{\partial}{\partial \theta} \] 
of smooth vector fields on $\bS^1$. 
Let $\bd := \partial_\theta$. 
We say that a unitary representation $(\pi, \cH)$ of $G$ is a {\it positive energy 
representation} if $-i\dd\pi(\bd)  \geq 0$. 
A unitary representations $(\pi,\cH)$ of $G$ is semibounded if and only if 
either $\pi$ or $\pi^*$ is of positive energy 
\cite[Thm.~8.3]{Ne10}, but unfortunately all these representations 
are trivial \cite[Thm.~8.7]{Ne10} (see also \cite{GO86} for the algebraic context). 

This is the same phenomenon that we already observed for hermitian groups 
 in Section~\ref{sec:8} (see also Remark~\ref{rem:2.2}(b)). 

Up to isomorphy, $\cV(\bS^1)$ has a unique central extension 
\[ \vir = \R c \oplus_\omega \cV(\bS^1), \qquad 
 \omega(f,g) = \int_{\bS^1} f'g'' - f'' g' \] 
called the {\it Virasoro algebra}. It integrates to a central Lie group extension 
\[ \1 \to \R \times \Z \to \Vir \to \Diff(\bS^1)_+ \to \1, \] 
called the {\it Virasoro\ group}. 
If $\hat G = \Vir$,  we write $\hat H \cong \R^2$ for the inverse image 
of the subgroup $H \cong \T$ of rigid rotations in $G$. 
Identifying $\bd$ with the corresponding element 
$(0,\bd)$ in $\vir$, we define positive energy representations as above and, 
using the classification of open invariant cones in $\vir$, 
one shows that semiboundedness is equivalent 
to $\pi$ or $\pi^*$ being of positive energy 
\cite[Thm~.8.15, Cor.~8.16]{Ne10}. 

To classify semibounded representations, we can not directly 
use the method of holomorphic induction 
in the form described in Section~\ref{sec:7} because 
$\Vir$ is neither an analytic 
Lie group nor locally exponential \cite{Ne06}. 
However, the manifold $\hat G/\hat H \cong G/H$ 
carries the structure of a complex manifold 
on which $G$ acts smoothly by holomorphic maps \cite{KY88}. 
In \cite{NS14} we show that this fact, together with some refinements 
due to L.~Lempert \cite{Le95} to make holomorphic induction work  for 
the passage from bounded representations of $\hat H$ to $\hat G =\Vir$. 
In particular, we show that the 
irreducible semibounded representations of positive energy 
are in one-to-one correspondence with the unitary highest weight 
representations, which have been classified in the 1980s 
(\cite{GKO86}, \cite{FQS86},  \cite{GW84,GW85}).

\section{Semiboundedness for solvable  groups} 
\mlabel{sec:11}

In Sections~\ref{sec:8} to \ref{sec:10}, we discussed semibounded representations 
for very specific classes of groups. The following theorem shows that there are 
structural obstructions for the existence of semibounded representations. 

\begin{theorem} {\rm(\cite{NZ13})} 
If $G$ is a connected Lie group and either nilpotent or 
$2$-step solvable (the commutator group is abelian), 
then all semibounded representations factor through the abelian group 
$G/(G,G)$. 
\end{theorem}

This theorem suggests to take a look at $3$-step solvable groups. 
Typical examples arise as follows. 
Let $(V,\omega)$ be locally convex symplectic vector space 
($\omega$ is assumed to be continuous and non-degenerate). 
Then the Heisenberg group 
\[ \Heis(V,\omega) = \R \times V, \quad 
(z,v)(z',v') := (z + z' + \textstyle{\frac{1}{2}}\omega(v,v'), v + v') \] 
is a Lie group. We consider a homomorphism 
$\alpha \: \R \to \Sp(V,\omega)$ defining a smooth action of $\R$ on 
$V$ and write $D := \alpha'(0) \in \fsp(V,\omega)$ for its infinitesimal generator. 
Then the {\it (generalized) oscillator group} 
\[ G = \Heis(V,\omega) \rtimes_\alpha \R, \qquad 
(h,t)(h',t') = (h \alpha_t(h'),t+t'),\quad 
\alpha_t(z,v) = (z, \alpha_t(v)) \] 
is a $3$-step solvable Lie group. 

In \cite{Ze14}, C.~Zellner shows that the classification of semibounded 
representations of such groups can be completely reduced to the case where 
$V$ is the Fr\'echet space of smooth vectors for a unitary one-parameter group 
$U_t = e^{itH}$ with $H \geq 0$ and $\omega(v,w) = \Im\la v,w\ra$. 
Again, every semibounded representation 
$(\pi,\cH)$ either satisfies the positive energy condition 
$\inf \Spec(-i\dd\pi(H)) > - \infty$ or its dual does.

\begin{theorem}[Uniqueness Theorem, \cite{Ze14, Ze15}] 
If $\Spec(H) \subeq [\eps,\infty[$ for some $\eps > 0$, then all irreducible semibounded 
positive energy representations of $G$ are Fock representations 
(classified by the lowest eigenvalue of $-i\dd\pi(H)$)
and every positive energy representation is type I and a direct integral 
of Fock representations. 
\end{theorem}

The situation becomes more complicated if $\inf (\Spec(H))= 0$. 
Then the semibounded representation theory of $G$ is no longer type I 
\cite{Ze14}. More precisely,
if $V$ is countably dimensional, then all representations of $\Heis(V,\omega)$ 
extend to semibounded representations of some oscillator group 
$G$, so that the classification problem (for all oscillator groups) 
is equivalent to the classification of the irreducible representations of the 
Canonical Commutation Relations (CCR), which is a ``wild'' problem. 

\section{Conclusion and perspectives} 

We conclude this article with a brief discussion of 
further developments and of directions in which the 
theory of semibounded representations can move from here. 

{\bf Non-type I representations:} We have seen above that the semibounded representations 
form a sector of the unitary representation theory 
of infinite dimensional Lie groups for which the full machinery 
of $C^*$-algebras is available on the abstract level. 
Presently, the theory is rather well developed for many important 
classes of Lie groups for which the corresponding representation 
theory is of type~I. Beyond type I representations, 
we have seen how to reduce the bounded representation theory 
of gauge groups to the representation theory of UHF 
algebras. However, for generalized oscillator groups 
similar results are still missing (Section~\ref{sec:11}). 
Here, and in many other cases, one would like to see a class 
of $C^*$-algebras, generalizing UHF algebras, that can host 
bounded and semibounded representations which are not of type I. 
Natural candidates for such algebras show up in the 
context of positive energy representations 
of gauge groups and have to be explored further \cite{JN15}. 

{\bf Positive energy representations:} 
Another important problem is a better understanding of the relations between 
the positive energy condition 
$-i\dd\pi(\bd) \geq 0$ for a specific element $\bd \in \g$ and 
the semiboundedness condition. We have seen above in 
Sections~\ref{sec:8}-\ref{sec:10}, for a well-chosen $\bd$, 
semiboundedness of $\pi$ often becomes equivalent to either 
$\pi$ of $\pi^*$ being of positive energy. 
This phenomenon also appears in \cite{JN15}, but 
in \cite{Ne15} and the closely related \cite{MN15}, the situation 
is more complicated. To understand these issues, we need a better theory 
of open invariant cones in infinite dimensional Lie algebras; 
see \cite{Ne10} for some first steps in this direction. 

{\bf Holomorphic induction from unbounded representations:}  
In Section~\ref{sec:7} we have seen how to obtain semibounded representations 
of a Lie group $G$ by holomorphic induction from bounded representations 
of a subgroup $H$. This method is strong enough to cover a large variety of cases, 
such as the ones discussed in Sections~\ref{sec:8}-\ref{sec:10}. 
However, generalizations are needed for other classes of groups because there are 
natural situations, where the representation of $H$ on a subspace 
$V$ generated by its intersection with $(\cH^\infty)^{\fp^-}$ is not bounded. 
In this case one has to work with weaker structures on the holomorphic 
bundles by either using the dense subspace $V^\infty$ as fibers of the bundle 
or by relaxing the smoothness of the structure on the bundle. 

{\bf Unitary representations of Lie supergroups:} In this survey we did not 
touch  representations  of Lie supergroups, although they are closely 
related to semibounded representations. 
In the unitary representation theory of Lie supergroups, the most important 
class of representations are those for which the corresponding unitary representation 
of the even part $G$ with Lie algebra $\g_{\oline{0}}$ 
is semibounded because we always have 
$-i\dd\pi([x,x]) \geq 0$ for any odd element $x$ in the Lie superalgebra. 
This has the consequence that the method of smoothing operators 
and the theory of semibounded representations applies particularly nicely 
to supergroups. We refer to \cite{NS15} for the construction of full 
 ``host algebras'' for finite dimensional Lie supergroups. We expect similar 
constructions to work for large classes of infinite dimensional Lie supergroups 
as well. 

\vspace*{1in} {\bf Acknowledgements.} 

We are most greatful to Bas Janssens, Hadi Salmasian and Christoph Zellner 
for careful reading of this article and for various suggestions 
that improved the presentation.

\section{References}

\renewcommand{\refname}{}    
\vspace*{-20pt}              

\end{document}